
\magnification1200
\input amstex.tex
\documentstyle{amsppt}

\hsize=12.5cm
\vsize=18cm
\hoffset=1cm
\voffset=2cm

\footline={\hss{\vbox to 2cm{\vfil\hbox{\rm\folio}}}\hss}
\nopagenumbers
\def\DJ{\leavevmode\setbox0=\hbox{D}\kern0pt\rlap
{\kern.04em\raise.188\ht0\hbox{-}}D}

\def\txt#1{{\textstyle{#1}}}
\baselineskip=13pt
\def\hf{{\textstyle{1\over2}}}
\def\a{\alpha}\def\b{\beta}
\def\d{{\,\roman d}}
\def\e{\varepsilon}

\def\G{\Gamma}
\def\k{\kappa}
\def\s{\sigma}
\def\t{\theta}
\def\={\;=\;}

\def\zt{\zeta(\hf+it)}

\def\D{\Delta}

\def\z{\zeta}

 \def\t{\theta}
\def\hf{{\textstyle{1\over2}}}
\def\txt#1{{\textstyle{#1}}}

\font\tenmsb=msbm10
\font\sevenmsb=msbm7
\font\fivemsb=msbm5
\newfam\msbfam
\textfont\msbfam=\tenmsb
\scriptfont\msbfam=\sevenmsb
\scriptscriptfont\msbfam=\fivemsb
\def\Bbb#1{{\fam\msbfam #1}}

\def \NN {\Bbb N}

\font\ff=cmr8
\def\txt#1{{\textstyle{#1}}}
\baselineskip=13pt

\font\teneufm=eufm10
\font\seveneufm=eufm7
\font\fiveeufm=eufm5
\newfam\eufmfam
\textfont\eufmfam=\teneufm
\scriptfont\eufmfam=\seveneufm
\scriptscriptfont\eufmfam=\fiveeufm
\def\mathfrak#1{{\fam\eufmfam\relax#1}}

\font\tenmsb=msbm10
\font\sevenmsb=msbm7
\font\fivemsb=msbm5
\newfam\msbfam
     \textfont\msbfam=\tenmsb
      \scriptfont\msbfam=\sevenmsb
      \scriptscriptfont\msbfam=\fivemsb
\def\Bbb#1{{\fam\msbfam #1}}

\def \NN {\Bbb N}

  \def\rightheadline{{\hfil{\ff
  On the Riemann zeta-function and the divisor problem III}\hfil\tenrm\folio}}

  \def\leftheadline{{\tenrm\folio\hfil{\ff
   Aleksandar Ivi\'c }\hfil}}
  \def\emptyheadline{\hfil}
  \headline{\ifnum\pageno=1 \emptyheadline\else
  \ifodd\pageno \rightheadline \else \leftheadline\fi\fi}

\topmatter
\title
ON THE RIEMANN ZETA-FUNCTION AND THE DIVISOR PROBLEM III
\endtitle
\author   Aleksandar Ivi\'c  \endauthor
\address
Aleksandar Ivi\'c, Katedra Matematike RGF-a
Universiteta u Beogradu, \DJ u\v sina 7, 11000 Beograd, Serbia
\endaddress
\email {\tt
ivic\@rgf.bg.ac.yu,  aivic\@matf.bg.ac.yu} \endemail
\dedicatory
\medskip
Dedicated to Prof. Imre K\'atai on the occasion of his seventieth birthday
\enddedicatory
\abstract
{Let $\D(x)$ denote the error term in the Dirichlet
divisor problem, and $E(T)$ the error term in the asymptotic
formula for the mean square of $|\zt|$. If
$E^*(t) = E(t) - 2\pi\D^*(t/2\pi)$ with $\D^*(x) =
 -\D(x)  + 2\D(2x) - \hf\D(4x)$ and we set
 $\int_0^T E^*(t)\d t = 3\pi T/4 + R(T)$, then we obtain
$$
R(T) = O_\e(T^{593/912+\e}), \;\int_0^TR^4(t)\d t \ll_\e T^{3+\e},
$$
and
$$
\int_0^TR^2(t)\d t = T^2P_3(\log T) + O_\e(T^{11/6+\e}),
$$
where $P_3(y)$ is a cubic polynomial in $y$ with positive leading
coefficient. }
\endabstract
\endtopmatter

\document
\head
1. Introduction and statement of results
\endhead

This paper is the continuation of the author's works [5], [6], where the analogy
between the Riemann zeta-function $\z(s)$ and the divisor problem
was investigated.
As usual, let the error term in the classical Dirichlet divisor problem be
$$
\D(x) \;=\; \sum_{n\le x}d(n) - x(\log x + 2\gamma - 1),
\leqno(1.1)
$$
and
$$
E(T) \;=\;\int_0^T|\zt|^2\d t - T\left(\log\bigl({T\over2\pi}\bigr) + 2\gamma - 1
\right),\leqno(1.2)
$$
where $d(n)$ is the number of divisors of
$n$, $\z(s)$ is the Riemann zeta-function, and $ \gamma = -\G'(1) = 0.577215\ldots\,$
is Euler's constant. In view of F.V. Atkinson's classical explicit formula
for $E(T)$ (see [1] and [3, Chapter 15]) it was known long ago that
there are analogies between $\D(x)$ and $E(T)$. However,
instead of the error-term function $\D(x)$ it
is more exact to work with the modified
function $\D^*(x)$ (see  M. Jutila [7], [8] and T. Meurman [10]), where
$$
\D^*(x) := -\D(x)  + 2\D(2x) - \hf\D(4x)
= \hf\sum_{n\le4x}(-1)^nd(n) - x(\log x + 2\gamma - 1),
\leqno(1.3)
$$
which is a better analogue of $E(T)$ than $\D(x)$.
M. Jutila (op. cit.) investigated both the
local and global behaviour of the difference
$$
E^*(t) \;:=\; E(t) - 2\pi\D^*\bigl({t\over2\pi}\bigr),
$$
and in particular in [8] he proved that
$$
\int_0^T(E^*(t))^2\d t \ll T^{4/3}\log^3T.\leqno(1.4)
$$
In the first part of the author's work [5] the bound in (1.4) was complemented
with the new bound
$$
\int_0^T (E^*(t))^4\d t \;\ll_\e\; T^{16/9+\e};\leqno(1.5)
$$
neither (1.4) or (1.5) seem to imply each other.
Here and later $\e$ denotes positive constants which are arbitrarily
small, but are not necessarily the same ones at each occurrence,
while $a \ll_\e b$ (same as $a = O_\e(b))$ means that
the $\ll$--constant depends on $\e$.
In the second part of the same work (op. cit.) it was proved that
$$
\int_0^T |E^*(t)|^5\d t \;\ll_\e\; T^{2+\e},\leqno(1.6)
$$
and some further results on higher moments of $|E^*(t)|$ were obtained as well.
In [6] the author sharpened (1.4) to
$$
\int_0^T (E^*(t))^2\d t \;=\; T^{4/3}P_3(\log T) + O_\e(T^{7/6+\e}),\leqno(1.7)
$$
where $P_3(y)$ is a polynomial of degree three in $y$ with
positive leading coefficient, and all the coefficients may be evaluated
explicitly.

\bigskip
The aim of the present work is to investigate the integral of $E^*(t)$.
More precisely, we define the error-term function $R(T)$ by the relation
$$
\int_0^T E^*(t)\d t = {3\pi\over 4}T + R(T).\leqno(1.8)
$$
We have (see  [2], [4] for the first formula and [14] for the second one)
$$
\int_0^T E(t)\d t = \pi T + G(T),\quad \int_0^T \D(t)\d t = {T\over 4} + H(T),
\leqno(1.9)
$$
where both $G(T), H(T)$ are $O(T^{3/4})$ and also $\Omega_\pm(T^{3/4})$
(for $g(x) > 0\;(x>x_0)\,$ $f(x) = \Omega(g(x))$ means that $f(x) = o(g(x))$
does not hold as $x\to\infty$, $f(x) = \Omega_\pm(g(x)))$ means that there
are unbounded sequences $\{x_n\},\,\{y_n\},\,$ and constants $A,B>0$ such
that $f(x_n) > Ag(x_n)$ and $f(y_n) < -Bg(y_n)$). Since
$$
\int_0^T \D(at)\d t = {1\over a}\int_0^{aT}\D(x)\d x
\qquad(a>0,\;T>0)\leqno(1.10)
$$
holds, it is obvious from (1.3), (1.9) and (1.10) that ${3\pi\over 4}$ is the
``correct" constant in (1.8), and that trivially one has the bound
$R(T) = O(T^{3/4})$, so that the problem is to improve it. We shall prove
\bigskip
THEOREM 1. {\it We have}
$$
R(T) = O_\e(T^{593/912+\e}), \quad {593\over912} = 0.6502129\ldots\,.
\leqno(1.11)
$$

\bigskip
THEOREM 2. {\it We have
$$
\int_0^TR^2(t)\d t = T^2P_3(\log T) + O_\e(T^{11/6+\e}),\leqno(1.12)
$$
where $P_3(y)$ is a cubic polynomial in $y$ with positive leading
coefficient, whose all coefficients may be explicitly evaluated.}

\medskip
The asymptotic formula (1.12) bears resemblance to (1.7), and it is proved
by a similar technique. The exponents in the error terms are, in both cases,
less than the exponent of $T$ in the main term by 1/6. This comes from the
use of (2.9) of Lemma 2.5, and in both cases the exponent of the error
term is the limit of the method. From (1.7) one obtains that $E^*(T)
= \Omega(T^{1/6}(\log T)^{3/2})$, which shows that $E^*(T)$ cannot be
too small. Likewise, (1.7) yields the following

\medskip
{\bf Corollary.} We have
$$
R(T) \= \Omega\Bigl(T^{1/2}(\log T)^{3/2}\Bigr).\leqno(1.13)
$$

\bigskip
THEOREM 3. {\it We have}
$$
\int_0^TR^4(t)\d t \ll_\e T^{3+\e}.
\leqno(1.14)
$$

\bigskip
It is rather difficult to
ascertain the true maximum order of magnitude of $R(T)$, but the omega-result
(1.13) makes it reasonable to believe that maybe it is $T^{1/2+o(1)}$
$(T\to\infty)$. It also seems reasonable to conjecture that
$$
R(T) = O_\e(T^{1/2+\e})\leqno(1.15)
$$
holds. If (1.15) is true,
then from Lemma 3, taking $H = T^{1/4}$, it would follow that
$$
E^*(T) \;\ll_\e\;T^{1/4+\e}\leqno(1.16)
$$
or equivalently
$$
E(T) = 2\pi\D^*\Bigl({T\over2\pi}\Bigr) + O_\e(T^{1/4+\e}).\leqno(1.17)
$$
By [4, Theorem 1.2] and (1.17) we have
$$
\eqalign{&
|\z(\hf + iT)|^2 \ll \log T \int_{T-1}^{T+1}|\zt|^2\d t + 1\cr&
\ll \log T\Bigl(\log T + E(T+1) - E(T-1)\Bigr)\cr&
\ll_\e \log T\left(\log T + 2\pi\D^*\Bigl({T+1\over2\pi}\Bigr) -
2\pi\D^*\Bigl({T-1\over2\pi}\Bigr)\right) + T^{1/4+\e} \ll_\e T^{1/4+\e},\cr}
$$
since, from (1.3) and $d(n) \ll_\e n^\e$,
$$
\D^*(T+H) - \D^*(T) = O(H\log T) +\hf\sum_{4T<n\le4(T+H)}(-1)^nd(n) \ll_\e HT^\e
$$
holds for $1 \ll H \ll T$. Therefore the conjectural (1.17) implies
 the hitherto unproved bound
$$
\z(\hf + iT) \;\ll_\e\;T^{1/8+\e}.\leqno(1.18)
$$
This significance of (1.18) shows the strength of the conjecture (1.15),
and the importance of the estimation of $R(T)$ and its mean values.

\smallskip
Furthermore we note that if (1.17) is true, then $\t = \rho$, where
$$
\t = \inf\Bigl\{\; c> 0\;:\;E(T) = O(T^c)\;\Bigr\},\quad
\rho = \inf\Bigl\{\; d> 0\;:\;\D(T) = O(T^d)\;\Bigr\}.
$$
Namely as $\t\ge1/4$ and $\rho\ge1/4$ are known to hold (this follows e.g., from
mean square results, see [4]) $\t = \rho$ follows from (1.17) and $\rho = \sigma$,
proved recently by Lau--Tsang [9], where
$$
\s = \inf\Bigl\{\; s> 0\;:\;\D^*(T) = O(T^s)\;\Bigr\}.
$$
The reader is also referred to M. Jutila [7] for a discussion on some related
implications. In any case our unconditional
results on $R(T)$ show, as is to be expected, that there is a lot of cancellation
in the mean sense between $E(T)$ and $2\pi\D^*(T/(2\pi))$, or in other words that
the function $E^*(T)$ is on the average much smaller than either
$E(T)$ or $2\pi\D^*(T/(2\pi))$.

\head
2. The necessary lemmas
\endhead

\bigskip
In this section we shall state the lemmas which are necessary
for the proof of our theorems. The first one brings forth a formula
for $\int_0^T E(t)\d t$, which is closely related to
F.V.  Atkinson's classical explicit
formula for $E(T)$ (see [1] or e.g., Chapter 15 of [3] or Chapter 2 of [4]).
This is due to J.L. Hafner and the author [2] (see also Chapter 3 of [4]).

\bigskip
LEMMA 1. {\it We have}
$$
\eqalign{
\int_0^T E(t)\d t &= \pi T + {1\over2} \Bigl({2T\over\pi}\Bigr)^{3/4}
\sum_{n\le T}(-1)^nd(n)n^{-5/4}e_2(T,n)\sin f(T,n)\cr&
-2\sum_{n\le c_0T}d(n)n^{-1/2}\Bigl(\log {T\over2\pi n}\Bigr)^{-2}
\sin \left(T\log \Bigl( {T\over2\pi n}\Bigr) - T + {1\over4}\pi\right)
 \cr&+ O(T^{1/4}),\cr}\leqno(2.1)
$$
{\it where $\,c_0 = {1\over2\pi} + \hf - \sqrt{{1\over4}+{1\over2\pi}}\,,$
$\,{\roman{ar\,sinh}}\,x = \log(x + \sqrt{1+x^2}\,),$ and
for $\,1\le n \ll T$, }
$$
\eqalign{
e_2(T,n) &= \left(1+{\pi n\over T}\right)^{-1/4}\left\{\left({2T\over\pi n}\right)^{1/2}
{\roman {ar\,sinh}}\left({\pi n\over2T}\right)^{1/2}\right\}^{-1/2}\cr&
= 1 + b_1{n\over T} + b_2\Bigl({n\over T}\Bigr)^2 + \ldots\,,\cr
f(T,n) &= 2T{\roman {ar\,sinh}}\,\bigl(\sqrt{\pi n/(2T)}\,\bigr) + \sqrt{2\pi nT
+ \pi^2n^2} - {\txt{1\over4}}\pi\cr&
=  -\txt{1\over4}\pi + 2\sqrt{2\pi nT} +
a_3n^{3/2}T^{-1/2} + a_5n^{5/2}T^{-3/2} +
a_7n^{7/2}T^{-5/2} + \ldots\,.\cr}\leqno(2.2)
$$

\bigskip
We also need a formula for the integral of $\D^*(x)$. From a classical
result of G.F.  Vorono{\"\i}  [14] (this also easily follows from pp. 90-91
of [3]) we have
$$
\int_0^X\D(x)\d x = {X\over4} + {X^{3/4}\over2\sqrt{2}\pi^2}
\sum_{n=1}^\infty d(n)n^{-5/4}\sin(4\pi\sqrt{nX}- {\txt{1\over4}}\pi)
+ O(1).
$$
To relate the above integral to the one of $\D^*(x)$ we proceed as
on pp. 472-473 of [3], using (1.3) and (1.10). In this way we are
led to

\bigskip
LEMMA 2. {\it We have}
$$
\eqalign{
\int_0^T\D^*(t)\d t &=
 {T^{3/4}\over2\sqrt{2}\pi^2} \sum_{n\le T^{2}}
(-1)^nd(n)n^{-5/4}\sin(4\pi\sqrt{nT}- {\txt{1\over4}}\pi)\cr&
+ O(T^{1/4}).
\cr}\leqno(2.3)
$$

\bigskip
We need also a result which relates $E^*(T)$ to its integral over a short
interval. This is

\bigskip
LEMMA 3. {\it For $T^\e \le H \ll T$ we have, for some constant $C>0$,}
$$
\eqalign{
E^*(T) &\le {1\over H}\int_T^{T+H} E^*(t)\d t + CH\log T,\cr
E^*(T) &\ge {1\over H}\int_{T-H}^T E^*(t)\d t - CH\log T.\cr}\leqno(2.4)
$$

\bigskip
{\bf Proof}. From (1.2) we have,
for $0 \le u\ll T$,
$$\eqalign{\cr
0 \le &\int_T^{T+u}|\zt|^2\d t = (T+u)\Bigl(\log\bigl({T+u\over2\pi}\bigr)
+ 2\gamma-1\Bigr) \cr&-  T\Bigl(\log\bigl({T\over2\pi}\bigr)
+ 2\gamma-1\Bigr) + E(T+u) - E(T).\cr}
$$
By the mean-value theorem this implies
$$
E(T) \le E(T+u) + O(u\log T),
$$
giving by integration
$$
E(T) \le {1\over H}\int_T^{T+H} E(t)\d t + CH\log T\qquad(1 \ll H\ll T,\,C>0).
\leqno(2.5)
$$
From (1.3) we have ($T^\e \le H \ll T$)
$$
\D^*(T) - {1\over H}\int_T^{T+H} \D^*(t)\d t \ll H\log T + {1\over
H}\int_T^{T+H}\sum_{4T<n\le4t}d(n)\d t \ll H\log T
\leqno(2.6)
$$
on applying a result of P. Shiu [13] on the values of multiplicative functions
in short intervals. It follows that
$$
\D^*(T) = {1\over H}\int_T^{T+H} \D^*(t)\d t
+ O(H\log T)\qquad(T^\e \le H \ll T).
$$
Hence
$$\eqalign{
2\pi\D^*\Bigl({T\over2\pi}\Bigr) &= {2\pi\over
H}\int_{T/2\pi}^{{T/2\pi}+H} \D^*(x)\d x + O(H\log T)\cr& = {1\over
H}\int_T^{T+2\pi H}\D^*\Bigl({t\over2\pi}\Bigr)\d t + O(H\log T)\cr&
= {2\pi\over  H }\int_T^{T+H}\D^*\Bigl({t\over2\pi}\Bigr)\d t+
O(H\log T),\cr} \leqno(2.7)
$$
on replacing $H$ by $H/2\pi$ in the last step. On combining (2.5)
and (2.7) we obtain the first inequality in (2.4), and the second
one follows analogously.
\bigskip

LEMMA 4. {\it If $1 \ll K \ll T^{3/4}$, $c_1 \ne 0, c_3, \ldots, c_{2L-1}$
are real constants, $L\ge1$ is fixed, and
$$
F(T,n) = c_1(Tn)^{1/2} + c_3n^{3/2}T^{-12} + \cdots + c_{2L-1}n^{L-1/2}T^{3/2-L},
$$
then for $(\k,\,\lambda)$ an exponent pair we have}
$$
\sum_{K<k\le K'\le2K}(-1)^kd(k){\roman e}^{F(T,k)i} \;\ll\;
T^{\k/2}K^{(1+\lambda)/2}\log T.\leqno(2.8)
$$

\bigskip
{\bf Proof}. The factor $(-1)^k$ is innocuous, and in fact can
be omitted, as was done in Chapter 7 of [3].  It suffices thus to consider
$$
S := \sum_{k\le K}d(k){\roman e}^{F(T,k)i}
= 2\sum_{m\le\sqrt{K}}\sum_{n\le{K/m}}{\roman e}^{F(T,mn)i} - \sum_{m\le\sqrt{K}}\sum_{n\le\sqrt{K}}
{\roman e}^{F(T,mn)i},
$$
where the familiar hyperbola method was applied. The sums over $n$ are split into
$\ll \log T$ subsums over the ranges $N<n\le N'\le 2N$. In view of ($C_{a,b}\ne0$)
$$
{\partial^{a+b}F(T,mn)\over (\partial m)^a (\partial n)^b}
\;\sim\; C_{a,b}T^{1/2}m^{1/2-a}n^{1/2-b}\qquad(a,b = 0,1,2\ldots\,,mn \ll T^{3/4})
$$
it follows that (see Chapter 2 of [3] for the definition and properties of
exponent pairs)
$$
\sum_{N<n\le N'\le 2N}{\roman e}^{F(T,mn)i}\;\ll\;\left({mT\over N}\right)^{\k/2}
N^\lambda.
$$
Hence we obtain
$$
S \ll \log T\sum_{m\le\sqrt{K}}(mT)^{\k/2}\left\{\left({K\over m}\right)^{\lambda
- \k/2} + K^{(2\lambda-\k)/4}\right\} \ll
T^{\k/2}K^{(1+\lambda)/2}\log T.
$$

\bigskip
LEMMA 5 (cf. Lemma 3 of [6]). {\it For $a > -\hf$ a constant we have}
$$
\sum_{n\le x}d^2(n)n^a = x^{a+1}P_3(\log x;a) +
O_\e(x^{a+1/2+\e}),\leqno(2.9)
$$
{\it where $P_3(y;a)$ is a polynomial of degree three in $y$
whose coefficients depend on $a$,
and whose leading coefficient equals $1/(\pi^2(a+1))$. All the
coefficients of $P_3(y;a)$ may be explicitly evaluated}.

\
\bigskip
The last lemma is a new result of O. Robert--P. Sargos [12] which will be needed
in the proof of  Theorem 3. This is

\bigskip
LEMMA 6. {\it Let $k\ge 2$ be a fixed
integer and $\delta > 0$ be given.
Then the number of integers $n_1,n_2,n_3,n_4$ such that
$N < n_1,n_2,n_3,n_4 \le 2N$ and}
$$
|n_1^{1/k} + n_2^{1/k} - n_3^{1/k} - n_4^{1/k}| < \delta N^{1/k}
$$
{\it is, for any given $\e>0$,}
$$
\ll_\e N^\e(N^4\delta + N^2).\leqno(2.10)
$$

\bigskip
\head
3. Proof of Theorem 1
\endhead
\bigskip
From Lemma 1 and Lemma 2 we deduce that
$$
\eqalign{
&R(T) = O(T^{1/2}\log T)\, +\cr&
+ {1\over2} \Bigl({2T\over\pi}\Bigr)^{3/4}
\sum_{n\le T}(-1)^nd(n)n^{-5/4}
\left\{e_2(T,n)\sin f(T,n)-\sin(2\sqrt{2\pi nT}-\pi/4)
\right\}.\cr}
\leqno(3.1)
$$
The sum over $n$ is written as
$$
\sum_{n\le T^{3/4}}= \sum_{n\le T^{1/3}} + \sum_{T^{1/3}<n\le T}
= {\sum}_1 + {\sum}_2,
$$
say. In $\sum_1$ we use the asymptotic expansion (2.2) (actually
$a_3 = {1\over6}\sqrt{2\pi^3}$) to infer that
$$
\eqalign{{\sum}_1& = a_3\sum_{n\le T^{1/3}}(-1)^nd(n)n^{-5/4}
\left\{n^{3/2}T^{-1/2}\cos f(T,n) +
O(n^{5/2}T^{-3/2})\right\}\cr&
= a_3T^{-1/2}\sum_{n\le T^{1/3}}(-1)^nd(n)n^{1/4}\cos f(T,n) +
O(T^{-3/4})\cr& = a_3T^{-1/2}{\sum}_3 + O(T^{-3/4}\log T),
\cr}\leqno(3.2)
$$
say. The important thing is that in $\sum_3$ we have the
increasing function $n^{1/4}$, and in $\sum_2$  the decreasing
function $n^{-5/4}$, while the exponential factors (up to a
constant) will be the same. This implies that the contributions
both in $\sum_2$ and $\sum_3$ will be dominated by the
contribution of $n \asymp T^{1/3}$. Thus essentially in the
explicit formula (3.1) for $R(T)$ the ``critical" values of $n$
will be when $n\asymp T^{1/3}$. In the truncated formula ($1\ll
N\ll T$)
$$
\D^*(x) = {1\over\pi\sqrt{2}}x^{1\over4}
\sum_{n\le N}(-1)^nd(n)n^{-{3\over4}}
\cos(4\pi\sqrt{nx} - {\txt{1\over4}}\pi) +
O_\e(x^{{1\over2}+\e}N^{-{1\over2}}),
\leqno(3.3)
$$
and in the integrated analogue for $E(T)$ (see (2.4)), if we want bounds
of the type $\D^*(x)\ll_\e x^{s+\e}$ ($E(T)\ll_\e T^{d+\e}$)
with $s < 1/3$ (resp. $d< 1/3$) we have to take $N$ in (3.3)
with $N = x^{1-2s}$. This implies that $1-2s>1/3$, hence the
``critical" values for $n$ in (3.3) will be larger than $x^{1/3}$,
whereas in (3.2) they are of the order $T^{1/3}$. The ``closeness"
of $E(T)$ and $2\pi\D^*(T/(2\pi))$ is basically induced by this
phenomenon.

To continue with the estimation of $R(T)$,  we use the Taylor
expansion of $f(T,n)$ (see (2.2)) with $L$ terms, where $L$ is
chosen so large that the tails of the series will make a negligible
contribution (i.e., $O(1)$). The sums in $\sum_2$ and $\sum_3$ are
split into $O(\log T)$ subsums of the form (2.8), after the removal
of the monotonic coefficients $n^{1/4}$ and $n^{-5/4}$ by partial
summation. Applying Lemma 4 it follows that
$$
T^{3/4}{\sum}_1 \ll T^{1/4}T^{\k/2}\log^2T\cdot
T^{{1\over3}({1\over4} + {1\over2} + {\lambda\over2})}
= T^{{1+\k\over2} + {\lambda\over6}}\log^2T,
$$
and in a similar way one has the bound
$$
T^{3/4}{\sum}_2 \;\ll\; T^{{1+\k\over2} + {\lambda\over6}}\log^2T.
$$
Therefore we have the bound
$$
R(T)\;\ll\;T^{{1+\k\over2} + {\lambda\over6}}\log^2T + T^{1/2}\log T
\ll\;T^{{1+\k\over2} + {\lambda\over6}}\log^2T,\leqno(3.4)
$$
since $0\le\k\le\hf\le\lambda\le1$.
Already the trivial exponent pair $(0,1)$ gives the bound $R(T)
\ll T^{2/3}\log^2T$, which improves the bound $R(T) \ll T^{3/4}$
that was mentioned in Section 1. The exponent in (3.4) does not
exceed 2/3 if
$$
3\k + \lambda \le 1 \leqno(3.5)
$$
holds, but it is not likely that the exponent
$593/912=0.6502129\ldots\;$ in (1.11) of Theorem 1 can be
attained in this fashion. To attain this exponent (more
sophisticated present-day estimates can yield a slightly smaller
exponent), one has to use estimates for two-dimensional
exponential sums. In particular, for (1.11) one can use the bound
of G. Kolesnik, worked out in Chapter 7 of [3].  This is $(c\ne0,\;K\ll T^{1/2})$
$$
\sum_{K\le k\le K'\le2K}(-1)^kd(k){\roman e}^
{ic(kT)^{1/2}+idk^{3/2}T^{-1/2}} \ll_\e T^\e(T^{-{1\over16}}K^{{173\over152}}
+ T^{1\over16}K^{119\over152}),\leqno(3.6)
$$
as the terms  $a_5k^{5/2}T^{-3/2}+\ldots$ in the (2.2) make a
smaller contribution. The terms $n>T^{1/2}$ in (3.1) may be estimated by
Lemma 4 with $(\k,\lambda) = (2/18,13/18) = ABA(1/6,2/3)$ in the terminology
of exponent pairs. The contribution is seen to be $O_\e(T^{11/18+\e})$, $11/18
= 0.6111\ldots\,$.

In the bound (3.6) it is the first term on
the right-hand side that will make the larger contribution, which
is  found, similarly as in the derivation of (3.4), to be
$$\eqalign{&
\ll_\e T^{3/4+\e}\left\{\max_{K\le T^{1/3}}T^{-1/2-1/16}K^{1/4+173/152}
+ \max_{K\ge T^{1/3}}T^{-1/16}K^{-5/4+173/152}\right\} \cr&\ll_\e
T^{593/912+\e}.\cr}
$$

\head
4. Proof of Theorem 2
\endhead
Combining Lemma 1 and Lemma 2 we obtain
$$
\eqalign{
&R(T) =
{1\over2}\left({2T\over\pi}\right)^{3/4}
\sum_{T<n\le T^{2}}(-1)^{n+1}d(n)n^{-5/4}\sin(2\pi\sqrt{2nT}- {\txt{1\over4}}\pi)\cr&
+ {1\over2}\left({2T\over\pi}\right)^{3/4} \sum_{n\le T}
(-1)^nd(n)n^{-5/4}\left\{e_2(T,n)\sin f(T,n)-\sin(2\pi\sqrt{2nT}-
{\txt{1\over4}}\pi)\right\}\cr&
-2\sum_{n\le c_0T}d(n)n^{-1/2}\Bigl(\log {T\over2\pi n}\Bigr)^{-2}
\sin \left(T\log \Bigl( {T\over2\pi n}\Bigr) - T + {\txt{1\over4}}\pi\right)
+ O(T^{1/4}).
\cr}\leqno(4.1)
$$
We set, for $T\le t\le 2T$,
$$
\eqalign{
S_1(t) &:= \sum_{T<n\le T^{2}}(-1)^{n+1}d(n)n^{-5/4}
\sin(2\pi\sqrt{2nt}- {\txt{1\over4}}\pi)\cr
S_2(t) &:= 2\sum_{n\le c_0T}d(n)n^{-1/2}\Bigl(\log {t\over2\pi n}\Bigr)^{-2}
\sin \left(t\log \Bigl( {t\over2\pi n}\Bigr) - t + {\txt{1\over4}}\pi\right).\cr}
$$
\medskip
Note that the mean square bound ($c\ne0$)
$$
\eqalign{&
\int_T^{2T}\Bigl|\sum_{K< k\le2K}(-1)^kd(k){\roman e}^{\sqrt{ckt}i}\Bigr|^2\d t
\cr&
= T\sum_{K< k\le2K}d^2(k) + \sum_{K< m\ne n\le2K}(-1)^{m+n}d(m)d(n)
\int_T^{2T}{\roman e}^{\sqrt{ct}(\sqrt{m}-\sqrt{n})i}\d t\cr&
\ll TK\log^3T + \sqrt{T}\sum_{K< m\ne n\le2K}{d(m)d(n)\over|\sqrt{m}-\sqrt{n}|}\cr&
\ll_\e TK\log^3T + T^{1/2+\e}\sum_{K< m\ne n\le2K}{K^{1/2}\over|m-n|}\cr&
\ll_\e T^\e(TK + T^{1/2}K^{3/2}) \cr}\leqno(4.2)
$$
holds for $1\ll K \ll T^C\;(C>0)$, where we used the first derivative test
(see Lemma 2.1 of [3]). The same bound also holds if in the exponential we have
$f(t,k)$ (cf. (2.2)) instead of $\sqrt{ctk}$, as shown e.g.,
in the derivation of the mean square
formula for $E(t)$ in Chapter 15 of [3]. Using (4.2) it follows that
$$
\int_T^{2T} t^{3/2}S_1^2(t)\d t \ll_\e T^{1+\e},\leqno(4.3)
$$
and, similarly as in Chapter 15 of [3], one obtains
$$
\int_T^{2T} S_2^2(t)\d t \ll_\e T^{1+\e}.\leqno(4.4)
$$
By using (4.2) we also have the crude bound
$$
\int_T^{2T} R^2(t)\d t \ll_\e T^{2+\e}.\leqno(4.5)
$$
Therefore from (4.1) and (4.3)--(4.5) we infer, by using the
Cauchy-Schwarz inequality for integrals, that
(setting $A={1\over\pi\sqrt{2\pi}}$ for brevity)
$$
\eqalign{&
\int_T^{2T} R^2(t)\d t =O_\e(T^{7/4+\e}) \,+\cr&
A\int\limits_T^{2T}t^{3/2}\left(\sum_{n\le T}
(-1)^nd(n)n^{-{5\over4}}\left\{e_2(T,n)\sin f(T,n)-\sin(2\pi\sqrt{2nT}-
{\txt{1\over4}}\pi)\right\}\right)^2\d t.\cr}
\leqno(4.6)
$$
Further, for a given $\delta>0$ we split
$$
\sum_{n\le T} = \sum_{n\le\delta\sqrt{T}} + \sum_{\delta\sqrt{T}<n\le T}
= S_3(t) + S_4(t),
$$
say. It follows from (4.6) that
$$
\eqalign{&
\int_T^{2T} R^2(t)\d t =
A\int\limits_T^{2T}t^{3/2}\Bigl(S_3^2(t) + S_4^2(t) + 2S_3(t)S_4(t)\Bigr)\d t
\cr& + O_\e(T^{7/4+\e}).\cr}
\leqno(4.7)
$$
Again, by (4.2), it is seen that the integral with $S_4^2(t)$ is absorbed
by the error term in (4.7).

\smallskip
Next, we consider the integral with $S_3(t)S_4(t)$ in (4.7), writing $m$ for the
integer variable in $S_3(t)$. If $m < T^{1/3}$, then we observe that
$$
\eqalign{
& \sin f(t,m) - \sin(2\sqrt{2\pi mt} - \pi/4) = \sum_{k=1}^\infty
{(y-y_0)^k\over k!}\sin(y_0 + \hf k\pi)\cr&
y = f(t,m),\; y_0 = 2\sqrt{2\pi mt} - \pi/4, \; y-y_0= d_3m^{3/2}t^{-1/2}
+ d_5m^{5/2}t^{-3/2} + \ldots\,.\cr}\leqno(4.8)
$$
Therefore  in
$$
\int_T^{2T}t^{3/2}\sum_{m<T^{1/3}}\ldots\sum_{\delta\sqrt{T}<n\le T}\ldots\d t
$$
we shall encounter  the exponential factor
$$
\exp\left(\pm i\left\{f(t,m)-\sqrt{8\pi mt}\right\}\right)
\exp\left(\pm i\left\{f(t,n)-\sqrt{8\pi nt}\right\}\right).\leqno(4.9)
$$
In the first exponential we use (4.8), and
the dominant contribution comes from the term $k=1$. The first derivative test
shows that the contribution  is
$$\eqalign{
&\ll T^2\sum_{m\le T^{1/3}}d(m)m^{-5/4}\cdot m^{3/2}T^{-1/2}
\sum_{\delta\sqrt{T}<n\le T}d(n)n^{-5/4}n^{-1/2}\cr& \ll
T^{3/2}T^{{5\over4}\cdot{1\over3}}T^{-{1\over2}\cdot{3\over4}}\log^2T
= T^{37/24}\log^2T.
\cr}
$$

\smallskip
In case when $T^{1/3} < m \le \delta\sqrt{T}$ in $S_3(t)$, we shall have
exponentials of the form
$$
\exp(\pm if(t,m)\pm i\sqrt{8\pi nt}\,),\quad \exp(\pm if(t,m)\pm if(t,n)),
$$
$$
\exp(\pm i\sqrt{8\pi mt}\pm i\sqrt{8\pi nt}\,),\quad
\exp(\pm i\sqrt{8\pi mt}\pm if(t,n)),
$$
with all possible combinations of signs. The most interesting case is that of
$$
\exp(iF(t,m,n)),\quad F(t,m,n) := \sqrt{8\pi mt} - f(t,n),
$$
when
$$
\eqalign{
{\d\over\d t}F(t,m,n) &= \sqrt{2\pi m\over t} - 2{\roman{ar\,sinh}}\sqrt{\pi n\over 2t}\cr&
= \sqrt{2\pi\over t}(\sqrt{m} - \sqrt{n}\,) + c_3n^{3/2}t^{-3/2} + c_5n^{5/2}t^{-5/2}
+ \ldots\,.\cr}
$$
Here we have $|\sqrt{m} - \sqrt{n}\,|\gg \sqrt{n}$ for $|n-m|\gg n$, namely for $m \ll n$.
In that case the contribution is clearly, by the first derivative test, $\ll_\e T^{7/4+\e}$.
If $m \gg n$, this means that
$$
\delta\sqrt{T}<n \ll \delta\sqrt{T}.\leqno(4.10)
$$
Then we have
$$
\left|\sqrt{2\pi\over t}\Bigl(\sqrt{m} - \sqrt{n}\,\Bigr)\right|
\;\gg\; \left({n\over t}\right)^{3/2}
$$
for $|m-n| \gg n^2/T$, which certainly holds in view of (4.10) if $\delta>0$ is
sufficiently small, since $|m-n|\ge 1$ when $m\ne n$.
The total contribution of such pairs $m,n$ is
$$
\ll T^2\sum_{m\le \sqrt{T}}d(m)m^{-5/4}\sum_{n\ne m,n\asymp\delta\sqrt{T}}
{d(n)\over|m-n|}n^{-5/4}n^{1/2} \;\ll\; T^{7/4}.
$$
In a similar fashion it is seen that all other cases make the same total contribution
which is $\ll T^{7/4}$. Thus we have
$$
\int_T^{2T} R^2(t)\d t = A\int_T^{2T}t^{3/2} S_3^2(t)\d t + O_\e(T^{7/4+\e}).
$$
In $S_3(t)$ we replace $e_2(t,n)$ by 1 (see (2.2)), making an error which is absorbed in the error term
above. Thus it is shown that
$$
\eqalign{&
\int_T^{2T} R^2(t)\d t = O_\e(T^{7/4+\e}) \,+\cr&
+ A\int_T^{2T}t^{3/2}\Bigl(\sum_{n\le\delta\sqrt{T}}
(-1)^nd(n)n^{-5/4}\Bigl\{\sin f(t,n) - \sin(\sqrt{8\pi nt}-
{\txt{1\over4}}\pi)\Bigr\}\Bigr)^2\d t\cr&
= A\int_T^{2T}t^{3/2}\sum_{n\le\delta\sqrt{T}}d^2(n)
n^{-5/2}\left\{\sin f(t,n) - \sin(\sqrt{8\pi nt}-
{\txt{1\over4}}\pi)\right\}^2\d t\cr&
\,+ O_\e(T^{7/4+\e}).\cr}\leqno(4.11)
$$
Namely when we square out the first sum above, then we encounter diagonal terms
($m = n$) which account for the main contribution.
There are also the non-diagonal terms ($m\ne n$),
which are estimated similarly as in the preceding
case, and which make a total contribution of $O_\e(T^{7/4+\e})$.

\medskip

At this point we invoke the elementary formula
$$
\eqalign{
(\sin \a - \sin\b)^2& = \sin^2\a + \sin^2\b - 2\sin\a\sin\b\cr&
= 1 - \hf(\cos2\a + \cos2\b) + \cos(\a+\b)-\cos(\a-\b)\cr}
$$
with
$$
\a = f(t,n),\quad \b= \sqrt{8\pi nt} - {\txt{1\over4}}\pi,
$$
and insert it in (4.11). To deal with the contribution of
$$
- \hf(\cos2\a + \cos2\b) + \cos(\a+\b)
$$
we split the sum over $n$ in (4.11) at $n = T^{\rho},\,0<\rho<\hf$.
Using $|\sin\a-\sin\b| \le |\a-\b|$ for $n < T^{\rho}$ and the
first derivative test for the remaining $n$ we obtain a contribution
which is
$$
\eqalign{&
\ll T^{5/2}\sum_{n < T^{\rho}}d^2(n)n^{1/2}T^{-1} + T^2\sum_{n\ge T^{\rho}}d^2(n)n^{-2}\cr&
\ll T^{{3\over2}+{3\over2}\rho}\log^3T + T^{2-\rho}\log^3T \ll T^{9/5}\log^3T\cr}
\leqno(4.12)
$$
with the choice $\rho = 1/5$. Using $1-\cos(\a-\b) = 2\sin^2(\hf(\a-\b))$ and
invoking the asymptotic expansion (2.2) for $f(T,n)$, we have altogether
$$
\int\limits_T^{2T}R^2(t)\d t = {\sqrt{2}\over\pi\sqrt{\pi}}
\sum_{T^{1/5}\le n\le\delta\sqrt{T}}{d^2(n)\over n^{5/2}}
\int\limits_T^{2T}t^{3\over2}\sin^2(\hf a_1n^{3/2}t^{-1/2})\d t
+ O(T^{15\over8}\log^3T).
\leqno(4.13)
$$
In the integral above we make the change of variable
$$
 \hf a_1n^{3/2}t^{-1/2} = y,\; \d t = -\hf a_1^2n^3y^{-3}\d y.
$$
We set  $z = (4(y/a_1)^2T)^{1/3}$ so that, after changing the order of integration
and summation,  the main term on the right-hand side of (4.13) becomes
$$\eqalign{&
{1\over\pi\sqrt{2\pi}}\sum_{T^{1/5}\le n\le\delta\sqrt{T}}d^2(n)n^{-5/2}
\int_{{1\over2}a_1 n^{3/2}(2T)^{-1/2}}^{{1\over2}a_1n^{3/2}T^{-1/2}}
(\hf a_1)^{3}n^{9/2}y^{-3}\sin^2y\cdot a_1^2n^3y^{-3}\d y\cr&
= {a_1^5\over8\pi\sqrt{2\pi}}\int_{2^{-3/2}a_1 T^{-1/5}}^{{1\over2}a_1\delta^{3/2}T^{1/4}}
\sum_{\max(T^{1/5},z)\le n\le \min(\delta\sqrt{T},2^{1/3}z)}d^2(n)n^5\cdot
{\sin^2y\over y^6}\d y.\cr}\leqno(4.14)
$$
The range of summation for $n$ is the interval $[z,\, 2^{1/3}z]$ if
$$
y\;\in\; J,\quad J := \left[{1\over2} a_1T^{-1/8},\,
{\delta^{3/2}\over2\sqrt{2}}a_1T^{1/4}\,\right].
$$
If we replace the interval of integration in the second integral in (4.14) by $J$,
then by using
$$
|\sin y| \;\le\;\min(1,|y|)\leqno(4.15)
$$
it follows that the error that is made is $\ll_\e T^{9/5+\e}$.

Now we use Lemma 3 ((2.9) with $a=5$) to obtain that the  integral
over $J$ equals, with suitable constants $d_j, e_j$,
$$
\eqalign{&
\int_{2^{-3/2}a_1T^{-1/5}}^{{1\over2} a_1\delta^{3/2}T^{1/4}}
\left(T^{2}y^{4}\sum_{j=0}^3d_j\log^j(y^{2/3}T^{1/3})
+ O_\e(T^{11/6+\e}y^{11/3})\right)\,{\sin^2y\over y^6}\d y\cr&
= T^{2}\int_{2^{-3/2}a_1T^{-1/5}}^{{1\over2} a_1\delta^{3/2}T^{1/4}}
\left(\sum_{j=0}^3e_j\log^j(y^2T)\right)\,{\sin^2y\over y^{2}}\d y
+ O_\e(T^{11/6+\e}),\cr}\leqno(4.16)
$$
since by using (4.15) we have
$$
\eqalign{&
\int_{2^{-3/2}a_1T^{-1/5}}^{{1\over2}a_1\delta^{3/2}T^{1/4}}
T^{{11\over6}+\e}\cdot {\sin^2y\over y^{7/3}}\d y \cr&=
T^{{11\over6}+\e}\left(\int_{2^{-3/2}a_1T^{-1/5}}^1
{\sin^2y\over y^{7/3}} \d y +
\int_1^{{1\over2}a_1\delta^{3/2}T^{1/4}}
{\sin^2y\over y^{7/3}}\d y\right)\cr&
\ll T^{{11\over6}+\e}\left(\int_0^1 y^{-1/3}\d y + \int_1^\infty y^{-7/3}\d y\right)
\ll T^{{11\over6}+\e},\cr}
$$
which accounts for the error term in Theorem 2.
Replacing the segment of integration in the integral on the right-hand
side of (4.16) by $(0,\,\infty)$, we make an error
which is $\ll_\e T^{9/5+\e}$. Namely, for
$0 < \a < 1 < \b,\; j = 0,1,\ldots\;$ fixed, we have
$$
\int_\a^\b\,{\sin^2y\over y^{2}}\log^j\d y =
\int_0^\infty\,{\sin^2y\over y^{2}}\log^j\d y
+ O(\a) + O(\b^{-1}\log^j\b),\leqno(4.17)
$$
where we used again (4.15).
Hence the expression in (4.16) becomes, on using (4.17) with
$\a =  2^{-3/2}a_1T^{-1/5}$,
$\b = 2a_1\delta^{3/2}T^{1/4}$,
$$
T^{2}\sum_{j=0}^3b_j\log^jT +  O_\e(T^{11/6+\e})
$$
with some constants $b_j\,(b_3>0)$ which may be explicitly evaluated, and
Theorem 2 follows.

\head
5. Proof of Theorem 3
\endhead
\bigskip
The proof of the fourth moment estimate in (1.14) follows by employing the method of
[5] used in the proof of (1.5) and (1.6), and therefore we shall be brief.
The chief ingredient of the proof is Lemma 6
with $k=2$, since raising the sum $|\sum_{K< k\le2K}(-1)^kd(k)
{\roman e}^{\sqrt{ckt}i}|$ to the fourth power leads to expressions of the form
$\sqrt{n_1} + \sqrt{n_2} - \sqrt{n_3} - \sqrt{n_4}\;(n_j\in\NN)$
in the exponential. Care has also to be taken
when one takes the first two terms in the asymptotic expansion
(2.2) of $f(t,k)$, namely
 of the term $k^{3/2}t^{-1/2}$. This is achieved by using the approach of
M. Jutila [7] part II, as embodied in e.g., Lemma 5 of [5], part I. As already
explained, the major contribution will come from the terms $n\asymp T^{1/3}$
in (3.1). The contribution of the terms $n\le T^{1/3}$, corresponding to $\sum_1$ in
the proof of Theorem 1, will be
$$
\ll T\log T\max_{K\ll T^{1/3}}\int_T^{2T}\Bigl|\sum_{K<n\le K'\le2K}(-1)^nd(n)n^{1/4}
{\roman e}^{\sqrt{ckt}i}\Bigr|^4\d t,\leqno(5.1)
$$
and the integral is, up to a small error term,
$$\eqalign{&
\ll \max_{K\ll T^{1/3}} K\int_{T/2}^{5T/2}
\Bigl|{\mathop{\sum\nolimits^*}\limits_{K<m,n,k,l\le K'\le 2K}}\times\cr&\times
(-1)^{m+n+k+l}d(m)d(n)d(k)d(l)\exp(i\D\sqrt{t}\,)\Bigr|\d t,\cr}
\leqno(5.2)
$$
where $\sum^* $ means that $|\D| \le T^{\e-1/2}$ holds, and
$$
\D \;:=\; \sqrt{8\pi}(\sqrt{m} + \sqrt{n} - \sqrt{k} - \sqrt{l}\,).
$$
Now we use Lemma 6 ((2.10) with $k=2$, $\delta \asymp
K^{-1/2}|\D|$), estimating the integral on the right-hand side of
(5.2) trivially. It follows that the total contribution will be
$$
\eqalign{&
\ll_\e T^{1+\e}\max_{K\ll T^{1/3}}T(K^{9/2}T^{-1/2} + K^3)\cr&
\ll_\e T^{3/2+3/2+\e} + T^{3+\e} \ll_\e T^{3+\e},\cr}
$$
and the same final bound follows for
the contribution of the terms $n$ in (5.1) satisfying $n>T^{1/3}$. The proof
of Theorem 3 is complete.

\medskip

\vfill
\eject
\topglue1cm
\bigskip
\Refs
\bigskip

\item{[1]} F.V. Atkinson, The mean value of the Riemann zeta-function,
Acta Math. {\bf81}(1949), 353-376.

\item{[2]} J.L. Hafner and A. Ivi\'c, On the mean square of the Riemann
zeta-function on the critical line, J. Number Theory {\bf31}(1989), 151-191.

\item{[3]} A. Ivi\'c, The Riemann zeta-function, John Wiley \&
Sons, New York, 1985 (2nd ed. Dover, Mineola, New York, 2003).

\item{[4]} A. Ivi\'c, The mean values of the Riemann zeta-function,
LNs {\bf 82}, Tata Inst. of Fundamental Research, Bombay (distr. by
Springer Verlag, Berlin etc.), 1991.

\item{[5]} A. Ivi\'c, On the Riemann zeta-function and the divisor problem,
Central European J. Math. {\bf(2)(4)} (2004), 1-15, and II, ibid.
{\bf(3)(2)} (2005), 203-214.

\item{[6]} A. Ivi\'c, On the mean square of the zeta-function and
the divisor problem, Annales  Acad. Scien. Fennicae Mathematica, in press.

\item{[7]} M. Jutila, Riemann's zeta-function and the divisor problem,
Arkiv Mat. {\bf21}(1983), 75-96 and II, ibid. {\bf31}(1993), 61-70.

\item{[8]} M. Jutila, On a formula of Atkinson, in ``Coll. Math. Sci.
J\'anos Bolyai 34, Topics in classical Number Theory, Budapest 1981",
North-Holland, Amsterdam, 1984, pp. 807-823.

\item{[9]} Y.-K. Lau and K.-M. Tsang, Omega result for the mean square of the
Riemann zeta-function, Manuscripta Math. {\bf117}(2005), 373-381.

\item{[10]} T. Meurman, A generalization of Atkinson's formula to $L$-functions,
Acta Arith. {\bf47}(1986), 351-370.

\item{[11]} T. Meurman, On the mean square of the Riemann zeta-function,
Quart. J. Math. (Oxford) {\bf(2)38}(1987), 337-343.

\item{[12]} O. Robert and P. Sargos, Three-dimensional
exponential sums with monomials, J. reine angew. Math. {\bf591}(2006), 1-20.

\item{[13]} P. Shiu, A Brun--Titchmarsh theorem for multiplicative
functions, J. reine angew. Math. {\bf31}(1980), 161-170.

\item{[14]} G.F. Vorono{\"\i}, Sur une fonction transcendante et ses applications
\`a la sommation de quelques s\'eries, Ann. \'Ecole Normale (3){\bf21}(1904),
2-7-267 and ibid. 459-533.

\endRefs
\vskip1cm

\enddocument

\bye